%Plain TEX file.  Run throught TEX twice so that it puts the numbers in.
%Plain TEX file.  Delete before this line
%seventy-two%seventy-two%seventy-two%seventy-two%seventy-two%seventy-two
\font\largebf=cmbx10 scaled\magstep2
\font\tenmsa=msam10
\font\tenmsb=msbm10
\def\all{\hbox{for all}}
\def\and{\hbox{and}}
\def\bra#1#2{\langle#1,#2\rangle}

\def\bigcupn{\bigcup\nolimits}
\def\cite#1\endcite{[#1]}
\def\dom{\hbox{\rm dom}}
\def\ds{\displaystyle}
\def\dst{^{**}}
\def\D{{\cal D}}

\def\exs{\hbox{there exists}}

\def\f#1#2{{#1 \over #2}}
\def\G{{\cal G}}
\def\infn{\inf\nolimits}

\def\I{\hbox{\rm I \hskip - .5em I}}
\def\kbar{\overline k}
\def\lr{\Longrightarrow}
\def\L{{\cal L}}
\def\minn{\min\nolimits}
\def\M{{\cal M}}

\def\PCLSC{{\cal PCLSC}}
\def\qed{\hfill\hbox{\tenmsa\char03}}
\def\qlr{\quad\lr\quad}

\def\quand{\quad\and\quad}
\def\r{\hbox{\tenmsb R}}
\def\rbar{\,]{-}\infty,\infty]}

\def\R{{\cal R}}
\def\Sbar{\overline{S}}
\long\def\slant#1\endslant{{\sl#1}}
\def\st{\hbox{such that}}
\def\Stilde{\wt{S}}
\def\supn{\sup\nolimits}
\def\T{\bf T}

\def\toto{\ {\mathop{\hbox{\tenmsa\char19}}}\ }
\def\ts{\textstyle}
\def\Ttilde{\wt{T}}
\def\wh{\widehat}
\def\wt{\widetilde}
\def\({\big(}
\def\){\big)}
\def\[{\big[}
\def\]{\big]}
\def\defSection#1{}
\def\defCorollary#1{}
\def\defDefinition#1{}
\def\defExample#1{}
\def\defLemma#1{}
\def\defNotation#1{}
\def\defProblem#1{}
\def\defRemark#1{}
\def\defTheorem#1{}
\def\locno#1{}
\def\meqno#1{\eqno(#1)}
\def\nmbr#1{}
\def\Proof{\medbreak\noindent{\bf Proof.}\enspace}
\def\Proo{{\bf Proof.}\enspace}
\def\Signoff{}
%%%%%%%%%%%%%%%%%%%%%%%%%%%%%%%%%%%%%%%%%%%%%%%%%%%%%%%%%%%%%%%%%%%%%%%
\def \INTsec{0}
\def \PRELIMsec{1}
\def \RSlem{1}
\def \HBLlem{2}
\def \CONTsec{2}
\def \CONTone{1}
\def \CONTrem{3}
\def \CONEthm{4}
\def \CONEone{2}
\def \CONEtwo{3}
\def \CONEthree{4}
\def \CONEfour{5}
\def \CONEfive{6}
\def \CONErem{5}
\def \EXthm{6}
\def \EXone{7}
\def \EXtwo{8}
\def \BANsec{3}
\def \FATnot{7}
\def \EXcor{8}
\def \EXSTone{9}
\def \EXSTtwo{10}
\def \SZthm{9}
\def \SUMABone{11}
\def \SUMABtwo{12}
\def \SZcor{10}
\def \EEcor{11}
\def \REPFNsec{4}
\def \Qone{13}
\def \Qtwo{14}
\def \Qthree{15}
\def \BRlem{12}
\def \KBARlem{13}
\def \VZthm{14}
\def \VZone{16}
\def \VZtwo{17}
\def \VZthree{18}
\def \VZfour{19}
\def \VZfive{20}
\def \VZsix{21}
\def \FFATlem{15}
\def \REPMFsec{5}
\def \MASVZthm{16}
\def \SUMthm{17}
\def \SUMcor{18}
\def \Dnot{19}
\def \NOTone{22}
\def \RESthm{20}
\def \PARthm{21}
\def \PARone{23}
\def \PARcor{22}
\def \QUALsec{6}
\def \TRANSlem{23}
\def \DIFFlem{24}
\def \DIFFone{24}
\def \SANDthm{25}
\def \SANDone{25}
\def \SANDtwo{26}
\def \SANDthree{27}
\def \SANDfour{28}
\def \ARENS{1}
\def \GOSSEZ{2}
\def \MASBR{3}
\def \MASUE{4}
\def \MASNI{5}
\def \MASNR{6}
\def \MASD{7}
\def \RCA{8}
\def \RS{9}
\def \PANDM{10}
\def \HBL{11}
\def \HBM{12}
\def \QUADARCHIV{13}
\def \SSDMON{14}
\def \SZNZ{15}
\def \VZJCA{16}
\def \ZBOOK{17}
%%%%%%%%%%%%%%%%%%%%%%%%%%%%%%%%%%%%%%%%%%%%%%%%%%%%%%%%%%%%%%%%%%%%%%%
\magnification 1200
\headline{\ifnum\folio<2
{\hfil{\largebf Quadrivariate existence theorems}\hfil}
\else\centerline{\bf Quadrivariate existence theorems and strong representability}\fi}
\centerline{\largebf and strong representability}
\medskip
\centerline{\bf Stephen Simons}
\medskip
\centerline{\bf Abstract}
\medskip
\noindent
In this paper, we give conditions under which we can compute the conjugate of a convex function on the product of two Fr\'echet spaces defined in terms of another convex function on the product of two (possibly different) Fr\'echet spaces.   We use this result to give simple proofs of some (both old and new) results for Banach spaces, and deduce some (both old and new) stability results for strongly representable multifunctions.   We take as our starting point a result on closed convex cones in the product of two Fr\'echet spaces.
%:  Section \INTsec
\defSection \INTsec
\bigbreak
\centerline{\bf\INTsec.\quad Introduction}
\medskip
\noindent
We start off by discussing Section \BANsec:
\par
\noindent
$\bullet$\enspace  If $X$ and $W$ are Banach spaces (all Banach spaces in this paper will be \slant real\endslant),\break $h$ is a proper, convex lower semicontinuous function on $W \times W^*$ and $C \in \L(X,W)$, let us define the convex function $k$ on $X \times X^*$ by $k(x,x^*) := \inf h(Cx \times (C^{\T})^{-1}x^*)$.   In\break \cite\MASNR, Proposition 2.2\endcite, Marques Alves and Svaiter gave  sufficient conditions for a certain formula for the ``restriction'' of $k^*$ to $X^* \times X$ to hold.
\par
\noindent
$\bullet$\enspace If $f$ and $g$ are proper, convex lower semicontinuous functions on the product of two (possible unrelated) Banach spaces, let $k$ be the inf--convolution of $f$ and $g$ with respect to the second variable.   In \cite\SZNZ, Theorem 4.2\endcite, Simons and Z\u{a}linescu gave sufficient conditions for $k^*$ to be the exact inf--convolution of $f^*$ and $g^*$ with respect to the first variable.
\par
\noindent
$\bullet$\enspace Finally, let $X$ and $Y$ be Banach spaces, $f$ and $g$ be proper, convex lower semicontinuous functions on $X \times X^*$ and $Y \times Y^*$, respectively, and $A \in \L(X,Y)$.  Define the function $k$ on $X \times X^*$ by $k(x,x^*) := \infn_{y^* \in Y^*}\[f(x,x^* - A^{\T}y^*) + g(Ax,y^*)\]$.   In \cite\VZJCA, Theorem 3.4\endcite, Voisei and Z\u{a}linescu (essentially) gave sufficient conditions for a certain formula for $k^*$ to hold.
\par
The initial purpose of this paper is to unify \cite\MASNR, Proposition 2.2\endcite, \cite\SZNZ, Theorem 4.2\endcite\ and \cite\VZJCA, Theorem 3.4\endcite.   We establish such a result in Theorem \EXthm, even for Fr\'echet spaces.   Of course, such a result must use four variables.   We show in Corollary \EXcor\ how Theorem \EXthm\ leads to a generalization of \cite\MASNR, Proposition 2.2\endcite, in Corollary \SZcor\  how Theorem \SZthm\ (a consequence of Theorem \EXthm) gives \cite\SZNZ, Theorem 4.2\endcite, and in Corollary \EEcor\ how Theorem \SZthm\ gives two results involving two Banach spaces and their duals, which we will apply in Section \REPFNsec.   Although it is possible to give a direct proof of Theorem \SZthm, \(see\break \cite\QUADARCHIV, Theorem 3\endcite\), we have deduced it here from the structurally much simpler Theorem \EXthm.    We note that, in Theorem \EXthm, the linear map $D$ goes from a space that $h$ uses to a space that $k$ uses, while the linear map $C$ goes from a space that $k$ uses to a space that $h$ uses.   
\par
In Section \REPFNsec, we define \slant strongly representative functions\endslant\ on the product of a Banach space and its dual.   We apply Corollary \EEcor\ in Theorem \VZthm\ and obtain sufficient ``sandwiched closed subspace'' conditions for various operations on two strongly representative functions to give rise to a third.   One of the  cases of Theorem \VZthm(a) is the result of \cite\VZJCA\endcite\ referred to above.   For reasons that we explain below, this analysis will be continued in Section \QUALsec. 
\par
In Section \REPMFsec, we define \slant strongly representable multifunctions\endslant\  from a Banach space into its dual.   In Theorem \SUMthm, Corollary \SUMcor, Theorem \PARthm\ and Corollary \PARcor, we give sufficient conditions for the generalized sum and generalized parallel sum of strongly representable multifunctions to be strongly representable.   The results on sums appear in \cite\VZJCA\endcite, but the results on parallel sums seem to be new.   The incentive for the study of the stability properties of strongly representable multifunctions stems from the knowledge that a multifunction is strongly representable $\iff$ it is maximally monotone ``of type (NI)'' $\iff$ it is maximally monotone ``of type (D)'' $\iff$ it is maximally monotone ``of type (ED)'' $\iff$ it is maximally monotone ``of type (WD)'' $\iff$ it is maximally monotone ``of dense type'' and that such a multifunction is ``of type (FP)'', ``of type (FPV)'', and ``strongly maximally monotone'', and has a very strong Br\o nsted--Rockafellar property.
\par
In Section \QUALsec, we will discuss some of the subtler features of the ``sandwiched closed subspace'' existence results initially discussed in Theorem \VZthm.   This kind of result was first introduced in \cite\SZNZ, Theorem 5.5\endcite\ in order to give sufficient conditions for the sum of maximally monotone multifunctions on a reflexive Banach space to be maximally monotone that were sufficiently flexible to subsume sufficient conditions previously obtained by a number of different authors.   We will illustrate the phenomenon under discussion in this paper with reference to Theorem \VZthm(a).   Write $H_1 := \ts\bigcupn_{\lambda > 0}\lambda\[A\(\D(\M f)\) - \D(\M g)\]$ and\break $H_2 := \ts\bigcupn_{\lambda > 0}\lambda\[A(\pi_E\,\dom\,f)- \pi_F\,\dom\,g\]$.   We know \slant a priori\endslant\ that $H_1$ is a (not necessarily convex) cone and $H_2$ is a convex cone.   Theorem \VZthm(a) implies that if there exists a closed linear subspace $H$ such that $H_1 \subset H \subset H_2$ then $k^*$ satisfies (\VZthree) and $\kbar$ is strongly representative.   What we prove in Theorem \SANDthm(a) is that, in addition, $H_1 = H_2$, and thus $H$ is uniquely determined.   So, in a sense, Theorem \SANDthm(a) belongs in Section \REPFNsec, but it uses Theorem \MASVZthm\ and the argument of Theorem \SUMthm\ from Section \REPMFsec.   We can make similar comments about Theorem \SANDthm(b), but this uses, in addition, Theorem \RESthm\ from Section \REPMFsec. 
\par
As we observed in our discussion of Section \BANsec\ above, the starting point for the later analysis is Theorem \EXthm\ in Section \CONTsec.   There are different ways of obtaining this result.   As explained in Remark \CONErem, we could start from the Attouch--Brezis theorem, as in \cite\MASNR\endcite, \cite\MASD\endcite, \cite\SZNZ\endcite\ and \cite\VZJCA\endcite.   We have chosen to take a different approach here, because some of our later results generalize the Attouch--Brezis theorem \(cf.\ \cite\SZNZ, Remark 4.3\endcite\).   We approach Theorem \EXthm\ through the result on \slant contragredient convex cones\endslant\ of Theorem \CONEthm.   Two comments about this result are in order: firstly, as we observe in Remark \CONErem, it actually suffices to consider closed linear subspaces rather than  closed convex cones --- this loss of generality does not result in a shortening of the analysis;  secondly, Theorem \CONEthm\ can also be deduced from \cite\ZBOOK, Theorem 2.8.6(v)\endcite\ --- we have included a direct proof here because it seems to be a fundamental and basic result.   As explained in Remark \CONTrem, the introduction of the notation of \slant contragredient\endslant\ was to avoid the proliferation of inverse maps, and because the notation ``$^*$'' is already overused.   (This term is taken from classical matrix theory, where the contragredient of a matrix is the inverse of its adjoint.)   It also has some notational advantages, which are explained in Remark \CONTrem.
\par
Finally, we make some comments about the results that appear in Section \PRELIMsec\ that we use in our derivation of Theorem \CONEthm.   We need two ingredients:  a convenient way of applying Baire's theorem, and a convenient way of assuring the existence of an appropriate continuous linear functional.   The first of these roles is satisfied by Lemma \RSlem, which first appeared in \cite\RS\endcite, and the second of these roles is satisfied by Lemma \HBLlem, which is an easy consequence of a result that first appeared in \cite\HBL\endcite, under the name of the ``Hahn--Banach--Lagrange theorem''.
\par
The author would like to thank Benar Svaiter for sending him a preprint of \cite\MASNI\endcite, Maicon Marques Alves for sending him a preprint of \cite\MASD\endcite, and Constantin Z\u alinescu for making him aware of \cite\MASBR\endcite, and sending him a preprint of \cite\VZJCA\endcite.
\par
Part of the research in this paper was done while the author was visiting \slant Centre de Recerca Matem\`atica\endslant\ in Barcelona, Spain.   The author would like to thank the \slant Centre\endslant\ for its hospitality.
%:  Section \PRELIMsec
\defSection \PRELIMsec
\bigbreak
\centerline{\bf \PRELIMsec.\quad Preliminary results}
\medskip
\noindent
If $E$ is a topological vector space, we write $\PCLSC(E)$ for the set of all proper, convex lower semicontinuous functions from $E$ into $\rbar$ and $E^*$ for the dual space of $E$.   If $F$ is also a topological vector space, we write $\L(E,F)$ for the set of continuous linear operators from $E$ into $F$.   Lemma \RSlem\ below was first proved in Rodrigues--Simons, \cite\RS, Lemma 1, pp.\ 1072--1073\endcite,  and was shown in \cite8\endcite\ to be a convenient way of applying Baire's theorem to a number of existence theorems in functional analysis.
%:  Lemma \RSlem
\defLemma \RSlem
\medbreak
\noindent
{\bf Lemma \RSlem.}\enspace\slant Let $P$ and $Q$ be complete metrizable topological vector spaces, $f \in \PCLSC(P)$ and $A \in \L(P,Q)$.   Suppose that $L := \ts\bigcupn_{\lambda > 0}\lambda A(\dom\,f)$ is a closed linear subspace of $Q$.   Then there exists $\gamma \ge 0$ such that $\big\{A(p)\colon\ p \in P,\ f(p) \le \gamma\big\}$ is a neighborhood of $0$ in L.\endslant
\medskip
Lemma \HBLlem\ is an immediate conseqence of the ``Hahn--Banach--Lagrange theorem'' of \cite\HBL, Theorem 2.9, p.\ 153\endcite\ or \cite\HBM, Theorem 1.11, p.\ 21\endcite\ (which is, in turn, a consequence of the Mazur--Orlicz theorem) and the fact that any linear functional dominated by a continuous seminorm is continuous.
%:  Lemma \HBLlem
\defLemma \HBLlem
\medbreak\noindent
{\bf Lemma \HBLlem.}\enspace\slant Let $S$ be a continuous seminorm on a topological vector space $Q$, $C$ be a nonempty convex subset of a vector space,  $f\colon\ C \to \rbar$ be proper and convex, and $A\colon\ C \to Q$ be affine.   Then there exists $q^* \in Q^*$ such that $\infn_C\[q^* \circ A + f\] = \infn_C\[S\circ A + f\]$.\endslant
%:  Section\CONTsec
\defSection\CONTsec
\bigbreak
\centerline{\bf\CONTsec.\quad The contragredient of a convex cone in a product space}
\medskip
\noindent
If $Q$ and $R$ are Fr\'echet spaces and $G$ is a convex cone in $Q \times R$ (always with vertex $0$) then we define the \slant contragredient\endslant, $G^\sim$, of $G$ (see Remark \CONTrem\ below for an explanation of this term) to be the convex cone
$$\big\{(q^*,r^*) \in Q^* \times R^*\colon\ \all\ (q,r) \in G,\ \bra{q}{q^*} \ge \bra{r}{r^*}\big\}.\meqno\CONTone$$
If $r \in R$, we write $G_r$ for the ``horizontal'' section $\big\{q \in Q\colon\ (q,r) \in G\big\}$.   Similarly, if $r^* \in R^*$, we write ${G^\sim}_{r^*}$ for the ``horizontal'' section $\big\{q^* \in Q^*\colon\ (q^*,r^*) \in G^\sim\big\}$.   We also write $G_R := \bigcup_{r \in R}G_r$, which is the projection of $G$ onto $Q$.   We shall show in Theorem \CONEthm\ how these concepts lead to a useful and easily statable result.  
%:  Remark \CONTrem
\defRemark \CONTrem
\medbreak
\noindent
{\bf Remark \CONTrem.}\enspace If $Q$ and $R$ are finite dimensional spaces and $G$ is a convex cone in $Q \times R$ then the \slant adjoint\endslant\ of $G$ was defined by Rockafellar in \cite\RCA, Section 39, pp.\ 413--423\endcite, extending a corresponding definition for linear subspaces of $Q \times R$ that was made by Arens in \cite\ARENS\endcite.    (The development in \cite\RCA\endcite\ was in terms of certain multifunctions, called \slant convex processes\endslant, from $Q$ into $R$ rather than subsets of $Q \times R$.)   In classical matrix theory, the \slant contragredient\endslant\ of a matrix is the inverse of its adjoint.   Our definition of the contragredient of a cone in (\CONTone) was derived from these precedents and the obvious interpretation of the \slant inverse\endslant\ of a subset of $Q \times R$.   We have made this definition to avoid the proliferation of inverse maps and because the notation ``$^*$'' is already overused.   This notation also leads to the formal similarity between the argument set for $k$ in (\CONEtwo) and the argument set for $k^*$ in (\CONEthree). 
\medskip
Theorem \CONEthm\ can also be deduced from \cite\ZBOOK, Theorem 2.8.6(v), p.\ 125--126\endcite, which is couched in the language of convex processes.
%:  Theorem \CONEthm
\defTheorem \CONEthm
\medbreak
\noindent
{\bf Theorem \CONEthm.}\enspace\slant Let $Q$ and $R$ be Fr\'echet spaces, $G$ be a closed convex cone in $Q \times R$, $h \in \PCLSC(Q)$ and
$$\ts L := \bigcupn_{\lambda > 0}\lambda\[\dom\,h - G_R\]\ \hbox{be a closed linear subspace of}\ Q.\meqno\CONEone$$
For all $r \in R$, let
$$k(r) := \inf h(G_r) > -\infty.\meqno\CONEtwo$$
Let $r^* \in R^*$ and $k^*(r^*) < \infty$.   Then
$$k^*(r^*) = \min h^*({G^\sim}_{r^*}).\meqno\CONEthree$$
\endslant
%:     Proof of Theorem \CONEthm 
\Proo Let $q^* \in {G^\sim}_{r^*}$ and $r \in R$.   If $q \in G_r$ then (\CONTone) implies that $\bra{r}{r^*} \le \bra{q}{q^*} \le h(q) + h^*(q^*)$.   Taking the infimum over $q \in G_r$ and using (\CONEtwo), $\bra{r}{r^*} \le k(r) + h^*(q^*)$.   It follows easily from this that $k^*(r^*) \le h^*(q^*)$.   So what we must prove is that there exists $q^* \in {G^\sim}_{r^*}$ such that $h^*(q^*) \le k^*(r^*)$.   It suffices for this that $q^* \in Q^*$ and 
$$p \in \dom\,h\ \and\ (q,r) \in G \qlr \bra{p}{q^*} - h(p) + \bra{r}{r^*} - \bra{q}{q^*} \le  k^*(r^*).\meqno\CONEfour$$ 
Define $A \in \L(Q \times Q \times R,Q)$ by $A(p,q,r) = p - q$ and $f \in \PCLSC(Q \times Q \times R)$ by $f(p,q,r) := h(p) + \I_G(q,r) - \bra{r}{r^*} + k^*(r^*)$, where $\I_G$ is the indicator function of $G$.   Then $A(\dom\,f) = \dom\,h - G_R$.   Thus (\CONEone) implies that Lemma \RSlem\ can be applied, from which there exist a continuous seminorm $T$ on $Q$ and $\gamma \ge 0$ such that
$$\big\{l \in L\colon\ T(l) < 1\big\} \subset \big\{p_1 - q_1\colon\ p_1 \in \dom\,h,\ (q_1,r_1) \in G,\ h(p_1) - \bra{r_1}{r^*} + k^*(r^*) \le \gamma\big\}.$$
Let $p \in \dom\,h$ and $(q,r) \in G$, so that $q - p = -(p - q) \in L$.   Let  $\mu > T(q - p)$.   Then $T\((q - p)/\mu\) < 1$, and so the above inclusion provides $p_1 \in \dom\,h$ and $(q_1,r_1) \in G$ such that
$$p_1 - q_1 = (q - p)/\mu \quand h(p_1) - \bra{r_1}{r^*} + k^*(r^*) \le \gamma.\meqno\CONEfive$$
Write $\alpha := \mu/(1 + \mu)$ and $\beta := 1/(1 + \mu)$: then $\alpha q_1  + \beta q= \alpha p_1 + \beta p$.   Since
$$(\alpha q_1  + \beta q,\alpha r_1  + \beta r) = \alpha(q_1,r_1) + \beta(q,r) \in G,$$
(\CONEtwo) implies that $k(\alpha r_1  + \beta r) \le h(\alpha q_1  + \beta q) = h(\alpha p_1  + \beta p) \le \alpha h(p_1) + \beta h(p)$.   Thus
$$\eqalign{0 
&\le k(\alpha r_1  + \beta r) - \bra{\alpha r_1  + \beta r}{r^*} + k^*(r^*) \le \alpha h(p_1) + \beta h(p) - \bra{\alpha r_1  + \beta r}{r^*} + k^*(r^*)\cr
&= \alpha\[h(p_1) - \bra{r_1}{r^*} + k^*(r^*)\]  + \beta\[h(p) - \bra{r}{r^*} + k^*(r^*)\]\cr
&\le \alpha\gamma  + \beta\[h(p) - \bra{r}{r^*} + k^*(r^*)\],}$$
where the final inequality follows from (\CONEfive).   If we now multiply by $1 + \mu$, we obtain $0 \le \mu\gamma  + h(p) - \bra{r}{r^*} + k^*(r^*)$ and, letting $\mu \to  T(q - p)$, we see that
$$0 \le \gamma T(q - p)  + h(p) - \bra{r}{r^*} + k^*(r^*).$$
From Lemma \HBLlem\ \(with $C = Q \times Q \times R$ and $S := \gamma T$\), we obtain $q^* \in Q^*$ such that, for all $p \in \dom\,h$ and $(q,r) \in G$, $\bra{q - p}{q^*}  + h(p) - \bra{r}{r^*} + k^*(r^*) \ge 0$.   This gives (\CONEfour), which completes the proof of Theorem \CONEthm.\qed
\medskip
In what follows, we write $\cdot^{\T}$ for the adjoint of the map $\cdot$.
%:  Remark \CONErem
\defRemark \CONErem
\medbreak
\noindent
{\bf Remark \CONErem.}\enspace The only application of Theorem \CONEthm\ in the rest of this paper is in the proof of Theorem \EXthm.   In order to obtain this application, it suffices to assume in Theorem \CONEthm\ that $G$ is a closed linear subspace rather than a closed convex cone.   However, this dilution of the assumptions does not result in a shortening of the analysis.
\par
Now in \cite\MASNR\endcite, \cite\MASD\endcite, \cite\SZNZ\endcite\ and \cite\VZJCA\endcite\ the starting point for the analysis was the (Banach space version of the) Attouch--Brezis theorem:
\smallskip
\noindent
\slant Let $P$ and $Q$ be nonzero Fr\'echet spaces, $g \in \PCLSC(P)$, $h \in \PCLSC(Q)$ and $A \in L(P,Q)$.   Suppose that $\ts\bigcupn_{\lambda > 0}\lambda\[\dom\,h - A(\dom\,g)\]$ is a closed linear subspace of $Q$.   Then, for all $p^* \in P^*$, $(g + h \circ A)^*(p^*) = \min_{q^* \in Q^*}\[g^*\(p^* - A^{\T}(q^*)\) + h^*(q^*)\]$.\endslant
\smallskip
\noindent
This Fr\'echet space version was first proved in Rodrigues--Simons, \cite\RS, Theorem 6, p.\ 1076\endcite, and much more general results were established in Z\u{a}linescu, \cite\ZBOOK, Theorem 2.8.3, p.\ 123--124\endcite.   While we have avoided using this result in our development since some of our later results imply the Attouch--Brezis theorem, we thought that it would be interesting for some readers to see a proof of Theorem \CONEthm\ that uses this instead of Lemmas \RSlem\ and \HBLlem.   Here is such a proof: Since $\dom\,\I_G = G$, (\CONEone) implies that\quad $\ts\bigcup_{\lambda > 0}\lambda\[\dom\,h - \pi_Q(\dom\,\I_G)\]$\quad is a closed linear subspace of $Q$, where $\pi_Q\colon\ Q \times R \to Q$ is defined by $\pi_Q(q,r) := q$.   Thus
$$\eqalign{k^*(r^*)
&=\supn_{r \in R}\[\bra{r}{r^*} - \inf h(G_r)\] = \supn_{(q,r) \in G}\[\bra{r}{r^*} - h(q)\]\cr
&= (\I_G + h \circ\pi_Q)^*(0,r^*)
= \minn_{q^* \in Q^*}\[{\I_G}^*\((0,r^*) - {\pi_Q}^{\T}(q^*)\) + h^*(q^*)\].}
$$
Now ${\I_G}^*\((0,r^*) - {\pi_Q}^{\T}(q^*)\) = \supn_{(q,r) \in G}\[\bra{r}{r^*} - \bra{q}{q^*}\]$.   It is easily seen that if this quantity is finite then $(q^*,r^*) \in G^\sim$, and then it has the value $0$.   Thus we obtain (\CONEthree), which gives us Theorem \CONEthm.   While this proof is shorter than that given using Lemmas \RSlem\ and \HBLlem, it is not shorter if one takes into account the overhead required to go from Lemmas \RSlem\ and \HBLlem\ to the Attouch--Brezis theorem.   We note \(cf.\ \cite\SZNZ, Remark 4.3, p.\ 10\endcite\) that Theorem \SZthm\ generalizes the Attouch--Brezis theorem for Banach spaces. 
\medskip
Theorem \EXthm, which first appeared (for Banach spaces) in \cite\QUADARCHIV, Theorem 21, pp.\ 12--13\endcite, will be applied in Corollary \EXcor\ and Theorem \SZthm.   The following chart should help the reader keep track of the various spaces and maps.
$$\matrix{&&X&{\ds C \atop \ds\longrightarrow}& W\cr
\rbar&{\ds k \atop \ds \longleftarrow}&\times&&\times&{\ds h \atop \ds \longrightarrow}&\rbar.\cr
&&U&{\ds D \atop \ds \longleftarrow}& T}$$
In Theorem \EXthm, and Corollary \EXcor, $\pi_W\colon\ W \times T \to W$ is defined by $\pi_W(w,t) := w$.
%:  Theorem \EXthm
\defTheorem \EXthm
\medbreak
\noindent
{\bf Theorem \EXthm.}\enspace\slant Let $X$, $W$, $U$ and $T$ be Fr\'echet spaces, $C \in \L(X,W)$, $D \in \L(T,U)$,\break $h \in \PCLSC(W \times T)$ and
$$\ts\bigcupn_{\lambda > 0}\lambda\[\pi_W\dom\,h - C(X)\]\ \hbox{be a closed linear subspace of}\ W.\meqno\EXone$$
For all $(x,u) \in X \times U$, let
$$k(x,u) := \inf h(Cx \times D^{-1}u) > -\infty.\meqno\EXtwo$$
Then, for all $(x^*,u^*) \in X^* \times U^*$ with $k^*(x^*,u^*) < \infty$,
$$k^*(x^*,u^*) = \min h^*\((C^{\T})^{-1}x^* \times D^{\T}u^*\).$$\endslant
%:    Proof of Theorem \EXthm
\Proo Let $Q := W \times T$ and $R := X \times U$, and $G$ be the closed convex cone
$$\big\{\((Cx,t),(x,Dt)\)\colon\ x \in X,\ T \in T\big\} \subset Q \times R.$$
Now, using the notation of Theorem \CONEthm, $G_R = C(X) \times T$, and so   
$$\dom\,h - G_R =\[\pi_W\dom\,h - C(X)\] \times T.$$ 
Thus $\bigcup_{\lambda > 0}\lambda\[\dom\,h - G_R\] = \bigcup_{\lambda > 0}\lambda\[\pi_W\dom\,h - C(X)\] \times T$, and (\EXone) implies that this is a closed linear subspace of $Q$.   Clearly $G_{(x,u)} = Cx \times D^{-1}u$, and so the function $k$ as defined in (\EXtwo) is identical with the function $k$ as defined in (\CONEtwo).   Since
$$\eqalign{
(w^*,t^*) &\in {G^\sim}_{(x^*,u^*)}\cr
&\iff \(x \in X\ \and\ t \in T \lr \bra{Cx}{w^*} + \bra{t}{t^*} \ge \bra{x}{x^*} + \bra{Dt}{u^*}\)\cr
&\iff \(x \in X\ \and\ t \in T \lr \bra{Cx}{w^*} \ge \bra{x}{x^*}\ \and\ \bra{t}{t^*} \ge \bra{Dt}{u^*}\)\cr
&\iff \(x \in X \lr \bra{Cx}{w^*} \ge \bra{x}{x^*}\)\ \and\ \(t \in T \lr \bra{t}{t^*} \ge \bra{Dt}{u^*}\)\cr
&\iff C^{\T}w^* =  x^*\ \and\ t^* = D^{\T}u^*
\iff (w^*,t^*) \in (C^{\T})^{-1}x^* \times D^{\T}u^*,}$$
the result follows from Theorem \CONEthm.\qed
%:  Section\BANsec
\defSection\BANsec
\bigbreak
\centerline{\bf\BANsec.\quad Applications to Banach spaces}
\medskip
\noindent
We now specialize to linear maps on Banach spaces, not because continued development in Fr\'echet spaces would be impossible, but simply because defining biduals is much more technical in the Fr\'echet space case.  Corollary \EXcor\ is a considerable strengthening of Marques Alves--Svaiter \cite \MASNR, Proposition 2.2, pp.\ 557, 562--563\endcite\ in that it provides a formula for $h^*$ rather than only for $h^@$.   It is very curious that the argument set for $h$ and the argument set for $h^@$ in Corollary \EXcor\ are identical. 
%:  Notation \FATnot
\defNotation \FATnot
\medbreak
\noindent
{\bf Notation \FATnot.}\enspace The following notation is taken from the theory of SSD spaces (see \cite\PANDM\endcite, \cite\HBM\endcite\ and \cite\SSDMON\endcite):  if $E$ is a Banach space and $f \in \PCLSC(E \times E^*)$, $f^@(x,x^*)$ stands for $f^*(x^*,\wh x)$ where, for all $x \in E$, $\wh x$ is the canonical image of $x$ in the bidual, $E\dst$.
%:  Corollary \EXcor
\defCorollary \EXcor
\medbreak
\noindent
{\bf Corollary \EXcor.}\enspace\slant Let $X$ and $W$ be Banach spaces, $C \in \L(X,W)$, $h \in \PCLSC(W \times W^*)$ and
$$\ts\bigcupn_{\lambda > 0}\lambda\[\pi_W\dom\,h - C(X)\]\ \hbox{be a closed linear subspace of}\ W.$$
For all $(x,x^*) \in X \times X^*$, let
$$k(x,x^*) := \inf h(Cx \times (C^{\T})^{-1}x^*) > -\infty.$$
Then, for all $(x^*,x\dst) \in X^* \times X\dst$ with $k^*(x^*,x\dst) < \infty$,
$$k^*(x^*,x\dst) = \min h^*\((C^{\T})^{-1}x^* \times C^{\T\T}x\dst\),\meqno\EXSTone$$
and, for all $(x,x^*) \in X \times X^*$ with $k^@(x,x^*) < \infty$, 
$$k^@(x,x^*) = \min h^@\(Cx \times (C^{\T})^{-1}x^*\).\meqno\EXSTtwo$$\endslant
%:     Proof of Corollary \EXcor
\Proo (\EXSTone) is immediate from Theorem \EXthm\ with $U := X^*$, $T := W^*$ and $D := C^{\T}$, and (\EXSTtwo) follows from (\EXSTone) and the observation that $C^{\T\T}\wh x = \wh{Cx}$.\qed
\medbreak
We now come to Theorem \SZthm, our quadrivariate version of the Attouch--Brezis theorem, which first appeared in \cite\QUADARCHIV, Theorem 3, pp.\ 2--4\endcite, and which will be applied in Corollary  \SZcor\ and Corollary \EEcor.   The following chart should help the reader keep track of the various spaces and maps.
$$\matrix{&&X&{\ds A \atop \ds\longrightarrow}& Y\cr
\rbar&{\ds f,k \atop \ds \longleftarrow}&\times&&\times&{\ds g \atop \ds \longrightarrow}&\rbar.\cr
&&U&{\ds B \atop \ds \longleftarrow}& V}$$
In Theorem \SZthm, $\pi_X\colon\ X \times U \to X$ is defined by $\pi_X(x,u) := x$, and $\pi_Y\colon\ Y \times V \to Y$ is defined by $\pi_Y(y,v) := y$.   We point out a subtle but important difference between Theorem \SZthm\ and the previous existence theorems --- the formula for $k^*(x^*,u^*)$ holds even if $k^*(x^*,u^*) = \infty$.   The reason for this that we can take (for instance) $y^* = 0$ and the proof of Theorem \CONEthm\ shows that $k^*(x^*,u^*) \le f^*(x^*,u^*) + g^*(0,B^{\T}u^*)$, and so both sides of this inequality have the same value, $\infty$.   On the other hand, in Theorem \EXthm\ it can easily happen that $(C^{\T})^{-1}x^* \times D^{\T}u^* = \emptyset$, and the existence of a minimizer is then problematical.
%:  Theorem \SZthm
\defTheorem \SZthm
\medbreak
\noindent
{\bf Theorem \SZthm.}\enspace\slant Let $X$, $Y$, $U$ and $V$ be Banach spaces, $A \in \L(X,Y)$, $B \in \L(V,U)$,\break $f \in \PCLSC(X \times U)$ and $g \in \PCLSC(Y \times V)$.   Let 
$$\ts\bigcupn_{\lambda > 0}\lambda\[A(\pi_X\,\dom\,f) - \pi_Y\,\dom\,g\]\ \hbox{be a closed linear subspace of}\ Y.\meqno\SUMABone$$
For all $(x,u) \in X \times U$, let
$$k(x,u) := \infn_{v \in V}\[f(x,u - Bv) + g(Ax,v)\] > -\infty.\meqno\SUMABtwo$$
Then, for all $(x^*,u^*) \in X^* \times U^*$,
$$k^*(x^*,u^*) = \minn_{y^* \in Y^*}\[f^*(x^* -  A^{\T}y^*,u^*) + g^*(y^*,B^{\T}u^*)\].$$\endslant
%:     Proof of Theorem \SZthm
\Proo We apply Theorem \EXthm\ with $W := X \times Y$, $T := U \times V$, $Cx := (x,Ax)$,\break $D(u,v) := u + Bv$, and $h\((x,y),(u,v)\) := f(x,u) + g(y,v)$.   It is then easy to verify that the function $k$ as defined in (\SUMABtwo) coincides with the function $k$ as defined in (\EXtwo) and $\pi_W\dom\,h = \pi_X\,\dom\,f \times \pi_Y\,\dom\,g$.   Since $\pi_W\dom\,h - C(X)$ is the inverse image of $A(\pi_X\,\dom\,f) - \pi_Y\,\dom\,g$ under the element of $\L(W,Y)$ defined by $(x,y) \mapsto Ax - y$, (\EXone) is a consequence of (\SUMABone).   We now obtain the required result since
$$\eqalignno{\min h^*\((C^{\T})^{-1}x^* \times D^{\T}u^*\)
&= \min \big\{h^*\(w^*,D^{\T}u^*\)\colon\ C^{\T}w^* = x^*\big\}\cr
&= \min \big\{h^*\((z^*,y^*),(u^*,B^{\T}u^*)\)\colon\ z^* + A^{\T}y^* = x^*\big\}\cr
&= \min \big\{f^*(z^*,u^*) + g^*(y^*,B^{\T}u^*)\colon\ z^* + A^{\T}y^* = x^*\big\}\cr
&= \minn_{y^* \in Y^*}\[f^*(x^* -  A^{\T}y^*,u^*) + g^*(y^*,B^{\T}u^*)\].&\qed}$$
\par
Corollary \SZcor\ first appeared in Simons--Z\u{a}linescu, \cite\SZNZ, Theorem 4.2, pp.\ 9--10\endcite.   In Corollary \SZcor, $\pi_E\colon\ E \times F \to E$ is defined by $\pi_E(x,y) := x$, and the result is immediate from Theorem \SZthm, with $X = Y = E$, $U = V = F$, and $A$ and $B$ identity maps.   
%:  Corollary \SZcor
\defCorollary \SZcor
\medbreak
\noindent
{\bf Corollary \SZcor.}\enspace\slant Let $E$ and $F$ be Banach spaces and $f,g \in \PCLSC(E \times F)$.   Suppose that $\bigcupn_{\lambda > 0}\lambda\[\pi_E\,\dom\,f - \pi_E\,\dom\,g\]$ is a closed linear subspace of $E$ and, for all $(x,y) \in E \times F$,
$$k(x,y) := \infn_{v \in F}\[f(x,y - v) + g(x,v)\] > -\infty.$$
Then, for all $(x^*,u^*) \in E^* \times F^*$,
$$k^*(x^*,u^*) = \minn_{y^* \in E^*}\[f^*(x^* -  y^*,u^*) + g^*(y^*,u^*)\].$$\endslant
\par
Corollary \EEcor\ contains two results involving two Banach spaces and their duals.   These results will be applied in Theorem \VZthm.   Corollary \EEcor\ first appeared in \cite\QUADARCHIV, Theorem 5,\ pp.\ 4--5\endcite.   Corollary \EEcor(a) is immediate from Theorem \SZthm\ with $X = E$, $Y = F$, $U = E^*$ and $V = F^*$, and Corollary \EEcor(b) is immediate from Theorem \SZthm\ with $X = E^*$, $Y = F^*$, $U = E$ and $V = F$, and changing the order of the arguments of $f$, $g$, $h$ and $k$.
\par
In the sequel, let $\pi_E\colon E \times E^* \to E$ be defined by $\pi_E(x,x^*):= x$, $\pi_F\colon F \times F^* \to F$ be defined by $\pi_F(y,y^*) := y$, $\pi_{E^*}\colon E \times E^* \to E^*$ be defined by $\pi_{E^*}(x,x^*) := x^*$, and $\pi_{F^*}\colon F \times F^* \to F^*$ be defined by $\pi_{F^*}(y,y^*) := y^*$.  
%:  Corollary \EEcor
\defCorollary \EEcor
\medbreak
\noindent
{\bf Corollary \EEcor.}\enspace\slant Suppose that $E$ and $F$ are Banach spaces, $f \in \PCLSC(E \times E^*)$ and $g \in \PCLSC(F \times F^*)$.   
\par
\noindent
{\rm(a)}\enspace Let $A \in \L(E,F)$, $B \in \L(F^*,E^*)$, $\bigcupn_{\lambda > 0}\lambda\[A(\pi_E\,\dom\,f) - \pi_F\,\dom\,g\]$ be a closed linear subspace of $F$ and, for all $(x,x^*) \in E \times E^*$,
$$k(x,x^*) := \infn_{y^* \in F^*}\[f(x,x^* - By^*) + g(Ax,y^*)\] > -\infty.$$
Then, for all $(x^*,x\dst) \in E^* \times E\dst$,
$$k^*(x^*,x\dst) = \minn_{y^* \in F^*}\[f^*(x^* -  A^{\T}y^*,x\dst) + g^*(y^*,B^{\T}x\dst)\].$$
{\rm(b)}\enspace Let $A \in \L(E^*,F^*)$, $B \in \L(F,E)$, $\bigcupn_{\lambda > 0}\lambda\[A(\pi_{E^*}\,\dom\,f)- \pi_{F^*}\,\dom\,g\]$ be a closed linear subspace of $F^*$ and, for all $(x,x^*) \in E \times E^*$,
$$k(x,x^*) := \infn_{y \in F}\[f(x - By,x^*) + g(y,Ax^*)\] > -\infty.$$
Then, for all $(x^*,x\dst) \in E^* \times E\dst$,
$$k^*(x^*,x\dst) = \minn_{y\dst \in F\dst}\[f^*(x^*,x\dst -  A^{\T}y\dst) + g^*(B^{\T}x^*,y\dst)\].$$\endslant\par
%
%:  Section \REPFNsec Representative and strongly representative functions
\defSection \REPFNsec
\medbreak
\centerline{\bf\REPFNsec.\quad Strongly representative functions}
\medskip
\noindent
We start of by recalling some facts from convex analysis.   Let $X$ be a Banach space and $k\colon X \to \rbar$ be proper and convex.   We write $\kbar$ for the (convex) function on $X$ such that the epigraph of $\kbar$ is the closure of the epigraph of $k$ in $X \times \r$.   Clearly, if $q\colon X \to \r$ is continuous and $k \ge q$ on $X$ then $\kbar \ge q$ on $X$.   It is also easy to see that $\kbar^* = k^*$.  
\par
In what follows, $\D(\cdot)$ stands for the domain of the multifunction $\cdot$, $\R(\cdot)$ stands for the range of the multifunction $\cdot$ and $\G(\cdot)$ stands for the graph of the multifunction $\cdot$.   Let $E$ be a Banach space.   If $f \in \PCLSC(E \times E^*)$ and,
$$(x,x^*) \in E \times E^* \qlr f(x,x^*) \ge \bra{x}{x^*}\meqno\Qone$$
then we define the multifunction $\M f \colon E \toto E^*$ by the rule
$$\G(\M f) = \{(x,x^*) \in E \times E^*\colon\ f(x,x^*) = \bra{x}{x^*}\}.\meqno\Qtwo$$
A function $f$ is said to be \slant strongly representative\endslant\ on $E \times E^*$ if $f \in \PCLSC(E \times E^*)$, (\Qone) is satisfied and, for all $(x^*,x\dst) \in E^* \times E\dst$,
$$(x^*,x\dst) \in E^* \times E\dst \qlr f^*(x^*,x\dst) \ge \bra{x^*}{x\dst}.\meqno\Qthree$$
In line with (\Qtwo), we then define the multifunction $\M f^*\colon E^* \toto E\dst$ by the rule
$$G(\M f^*) = \{(x^*,x\dst) \in E^* \times E\dst\colon\ f^*(x^*,x\dst) = \bra{x^*}{x\dst}\}.$$
In the main result of this section, Theorem \VZthm, we show how two strongly representative functions can be combined to give rise to a third, but we first give two preliminary results.
%:  Lemma \BRlem
\defLemma \BRlem
\medbreak
\noindent
{\bf Lemma \BRlem.}\enspace\slant Let $E$ be a Banach space and $f \in \PCLSC(E \times E^*)$ be strongly representative.   Then:
\par
\noindent
{\rm(a)}\enspace   If $\alpha, \beta > 0$ and $f(x,x^*) < \bra{x}{x^*} + \alpha\beta$ then there exists $(y,y^*) \in \G(\M f)$ such that
$\|y - x\| < \alpha$ and $\|y^* - x^*\| < \beta$.
\par
\noindent
{\rm(b)}\enspace $\D(\M f) \subset \pi_E\,\dom\,f \subset \overline{\D(\M f)} \quand \R(\M f) \subset \pi_{E^*}\,\dom\,f \subset \overline{\R(\M f)}$.  
\endslant
%:     Proof of Lemma \BRlem
\Proof (a) was established in \cite\MASBR, Theorem 3.4, pp.\ 700--701\endcite\ in the proof of Theorem \MASVZthm\ below, and (b) is immediate from (a) and the observation that $\bra{x}{x^*} \in \r$.   \(Much stronger results than (a) are now known --- see \cite\SSDMON, Theorem 9.9, pp. 254--255\endcite\).\qed
%:  Lemma \KBARlem
\defLemma \KBARlem
\medbreak
\noindent
{\bf Lemma \KBARlem.}\enspace\slant Let $E$ be a Banach space, $k\colon E \times E^* \to \rbar$ be proper and convex, for all $(x,x^*) \in E \times E^*$, $k(x,x^*) \ge \bra{x}{x^*}$ and, for all $(x^*,x\dst) \in E^* \times E\dst$, $k^*(x^*,x\dst) \ge \bra{x^*}{x\dst}$.   Then $\kbar$ is strongly representative.\endslant
%:     Proof of Lemma \KBARlem
\Proof Since the function $(x,x^*) \mapsto \bra{x}{x^*}$ is continuous, for all $(x,x^*) \in E \times E^*$, $\kbar(x,x^*) \ge \bra{x}{x^*}$.   The result follows since $\kbar^* = k^*$.\qed
\medskip
Theorem \VZthm\ appeared implicitly in \cite\QUADARCHIV, Theorem 20,\ pp.\ 12\endcite.   The case of Theorem \VZthm(a) with $H := \ts\bigcupn_{\lambda > 0}\lambda\[A\(\D(\M f)\) - \D(\M g)\]$  first appeared in Voisei--Z\u{a}linescu, \cite\VZJCA, Theorem 3.4, p.\ 1024\endcite, but we could also take $H = \ts\bigcupn_{\lambda > 0}\lambda\[A(\pi_E\,\dom\,f) - \pi_F\,\dom\,g\]$ or $H = \ts\bigcupn_{\lambda > 0}\lambda\[A\(\D(\M f)\) - \pi_F\,\dom\,g\]$.   See \cite\HBM, Remark 32.4,  p.\ 130\endcite\ for other possibilities, and the motivation for this kind of result.   Of course, there are similar possibilities for Theorem \VZthm(b).   The \slant statement\endslant\ of Theorem \VZthm\ is patterned after the ``sandwiched closed subspace theorem'' of \cite\HBM, Theorem 32.2,  p.\ 129\endcite\ or \cite\SZNZ, Theorem 5.5, p. 13\endcite\ (which were specific to Fitzpatrick functions), but the \slant proof\endslant\ here is much simpler by virtue of Lemma \BRlem(b) above.   Theorem \VZthm\ will be ``continued'' in Theorem \SANDthm.
%:  Theorem \VZthm
\defTheorem \VZthm
\medbreak
\noindent
{\bf Theorem \VZthm.}\enspace\slant Suppose that $E$ and $F$ are Banach spaces and that $f \in \PCLSC(E \times E^*)$ and $g \in \PCLSC(F \times F^*)$ are strongly representative.
\par
\noindent
{\rm(a)}\enspace Suppose that $A \in \L(E,F)$, there exists a closed linear subspace $H$ of $F$ such that
$$A\(\D(\M f)\) - \D(\M g) \subset H \subset \ts\bigcupn_{\lambda > 0}\lambda\[A(\pi_E\,\dom\,f)- \pi_F\,\dom\,g\],\meqno\VZone$$
and, for all $(x,x^*) \in E \times E^*$,
$$k(x,x^*) := \infn_{y^* \in F^*}\[f(x,x^* - A^{\T}y^*) + g(Ax,y^*)\].\meqno\VZtwo$$
Then, for all $(x^*,x\dst) \in E^* \times E\dst$,
$$k^*(x^*,x\dst) = \minn_{y^* \in F^*}\[f^*(x^* -  A^{\T}y^*,x\dst) + g^*(y^*,A^{\T\T}x\dst)\].\meqno\VZthree$$
Furthermore, $\kbar$ is a strongly representative function on $E \times E^*$.
\medbreak
\noindent
{\rm(b)}\enspace Suppose that $B \in \L(F,E)$, there exists a closed linear subspace $H$ of $F^*$ such that
$$B^{\T}\(\R(\M f)\) - \R(\M g) \subset H \subset \ts\bigcupn_{\lambda > 0}\lambda\[B^{\T}(\pi_{E^*}\,\dom\,f) - \pi_{F^*}\,\dom\,g\],\meqno\VZfour$$
and, for all $(x,x^*) \in E \times E^*$,
$$k(x,x^*) := \infn_{y \in F}\[f(x - By,x^*) + g(y,B^{\T}x^*)\].\meqno\VZfive$$
Then, for all $(x^*,x\dst) \in E^* \times E\dst$,
$$k^*(x^*,x\dst) = \minn_{y\dst \in F\dst}\[f^*(x^*,x\dst -  B^{\T\T}y\dst) + g^*(B^{\T}x^*,y\dst)\].\meqno\VZsix$$
Furthermore, $\kbar$ is a strongly representative function on $E \times E^*$.
\endslant
%:     Proof of Theorem \VZthm
\Proof(a)\enspace It is clear from Lemma \BRlem(b), the continuity of $A$ and (\VZone) that
$$A(\pi_E\,\dom\,f) - \pi_F\,\dom\,g \subset A\(\overline{\D(\M f)}\) - \overline{\D(\M g)} \subset \overline{A\(\D(\M f)\) - \D(\M g)} \subset H,$$
and so $\bigcupn_{\lambda > 0}\lambda\[A(\pi_E\,\dom\,f) - \pi_F\,\dom\,g\] = H$.   Since, for all $(x,x^*) \in E \times E^*$,
$$k(x,x^*) \ge \infn_{y^* \in F^*}\[\bra{x}{x^* - A^{\T}y^*} + \bra{Ax}{y^*}\] = \bra{x}{x^*} > -\infty,$$
the required formula for $k^*(x^*,x\dst)$ now follows from Corollary \EEcor(a) with $B := A^{\T}$.   Furthermore, for all $(x^*,x\dst) \in E^* \times E\dst$ and $y^* \in F^*$,
$$f^*(x^* -  A^{\T}y^*,x\dst) + g^*(y^*,A^{\T\T}x\dst) 
\ge \bra{x^* -  A^{\T}y^*}{x\dst} + \bra{y^*}{A^{\T\T}x\dst} = \bra{x^*}{x\dst}.$$
Consequently, $k^*(x^*,x\dst) \ge \bra{x^*}{x\dst}$.   From Lemma \KBARlem, $\kbar$ is strongly representative.
\par
(b)\enspace It is clear from Lemma \BRlem(b), the continuity of $B^{\T}$  and (\VZfour) that
$$B^{\T}(\pi_{E^*}\,\dom\,f) - \pi_{F^*}\,\dom\,g \subset B^{\T}\(\overline{\R(\M f)}\) - \overline{\R(\M g)} \subset \overline{B^{\T}\(\R(\M f)\) - \R(\M g)} \subset H,$$
and so $\ts\bigcupn_{\lambda > 0}\lambda\[B^{\T}(\pi_{E^*}\,\dom\,f) - \pi_{F^*}\,\dom\,g\] = H$.   Since, for all $(x,x^*) \in E \times E^*$,
$$k(x,x^*) \ge \infn_{y \in F}\[\bra{x - By}{x^*} + \bra{y}{B^{\T}x^*}\] = \bra{x}{x^*} > -\infty,$$
the required formula for $k^*(x^*,x\dst)$ now follows from Corollary \EEcor(b) with $A := B^{\T}$.   Furthermore, for all $(x^*,x\dst) \in E^* \times E\dst$ and $y\dst \in F\dst$,
$$f^*(x^*,x\dst -  B^{\T\T}y\dst) + g^*(B^{\T}x^*,y\dst)
\ge \bra{x^*}{x\dst -  B^{\T\T}y\dst} + \bra{B^{\T}x^*}{y\dst} = \bra{x^*}{x\dst}.$$
Consequently, $k^*(x^*,x\dst) \ge \bra{x^*}{x\dst}$.   From Lemma \KBARlem, $\kbar$ is strongly representative.\qed
\medskip
We end this section with another important computational result for strongly representative functions.   We refer the reader to Notation \FATnot\ for the definition of $f^@$.
%:  Lemma \FFATlem
\defLemma \FFATlem
\medbreak
\noindent
{\bf Lemma \FFATlem.}\enspace\slant Let $E$ be a Banach space and $f \in \PCLSC(E \times E^*)$ be strongly representative.   Then:
\par
\noindent
{\rm(a)}\enspace $f^@$ is strongly representative and $\M f^@ = \M f$.
\par
\noindent
{\rm(b)}\enspace $(s,s^*) \in \G(\M f) \iff (s^*,\wh s) \in \G(\M f^*)$. 
\endslant
%:     Proof of Lemma \FFATlem
\Proof (a) was established in  \cite\MASBR, Theorem 4.2(1), p.\ 702\endcite\ and \cite\VZJCA, Corollary 2.14, p.\ 1019\endcite.   \(A much more general result than (a) was proved in \cite\SSDMON, Theorem 5.8, p.\ 241--242\endcite.\)   (b) follows from (a) since
$$\eqalignno{(s,s^*) \in \G(\M f)
&\iff (s,s^*) \in \G(\M f^@) \iff f^@(s,s^*) = \bra{s}{s^*}\cr
&\iff f^*(s^*,\wh s) = \bra{s^*}{\wh s} \iff (s^*,\wh s) \in \G(\M f^*).&\qed}$$\par
\vfill
\eject
%:  Section \REPMFsec Strongly representable multifunctions
\defSection \REPMFsec
\par
\centerline{\bf\REPMFsec.\quad Strongly representable multifunctions}
\medskip
\noindent
A multifunction $S\colon\ E \toto E^*$ is said to be \slant strongly representable\endslant\ if there exists a strongly representative function $f \in \PCLSC(E \times E^*)$ such that $\M f = S$.   If $S$ is strongly representable, we say that $f$ is a \slant strong representor\endslant\ for $S$ if $f \in \PCLSC(E \times E^*)$ and $\M f = S$.   
\medskip
The following result was first proved in Marques Alves--Svaiter, \cite\MASBR, Theorem 4.2, pp.\ 702--704\endcite\ and Voisei--Z\u{a}linescu,\cite\VZJCA, Theorem 2.12\endcite.
%:  Theorem \MASVZthm
\defTheorem \MASVZthm
\medskip
\noindent
{\bf Theorem \MASVZthm.}\enspace \slant Any strongly representable multifunction is maximally monotone.\endslant
\medskip
After these two seminal papers, there were a number of improvements to Theorem \MASVZthm.   It was proved in Marques Alves--Svaiter, \cite\MASNI, Theorem 1.2, pp.\ 885, 887--889\endcite\ and\break \cite\MASD, Theorem 4.4, pp.\ 1084--1085\endcite\ that a multifunction is strongly representable $\iff$ it is maximally monotone ``of type (NI)'' $\iff$ it is maximally monotone ``of type (D)''.   It was then proved in \cite\SSDMON, Theorems 9.5, 9.7, p. 253\endcite\ that a multifunction is strongly representable $\iff$ it is maximally monotone ``of type (ED)'' $\iff$ it is maximally monotone ``of type (WD)'' $\iff$ it is maximally monotone ``of dense type''.   It is shown in \cite\SSDMON, Theorems 9.9 and 9.10, pp. 254--256\endcite\ that these results imply (among other properties) that such a multifunction is ``of type (FP)'', ``of type (FPV)'', and ``strongly maximally monotone'', and that they have a very strong Br\o nsted--Rockafellar property.   \(The first three of these results appeared in Voisei--Z\u{a}linescu, \cite\VZJCA, Remark 3.6, Theorems 4.1 and 4.2, p. 1024, pp.\ 1027--1030\endcite.\)   All these facts provide an incentive for the study of the stability properties of strongly representable multifunctions, which we will perform in Theorem \SUMthm, Corollary \SUMcor, Theorem \PARthm\ and Corollary \PARcor.   Theorem \SUMthm\ and Corollary \SUMcor\ appeared in Voisei--Z\u{a}linescu, \cite\VZJCA, Theorem 3.4 and Corollary 3.5, p. 1034\endcite.
%:  Theorem \SUMthm
\defTheorem \SUMthm
\medbreak
\noindent
{\bf Theorem \SUMthm.}\enspace\slant Let $E$ and $F$ be Banach spaces, $S\colon E \toto E^*$ and $T\colon F \toto F^*$ be strongly representable, $A \in \L(E,F)$ and $\ts\bigcupn_{\lambda > 0}\lambda\[A\(\D(S)\) - \D(T)\]$ be a closed linear subspace of F.  Then the multifunction $S + A^{\T}TA$ is strongly representable.\endslant
%:     Proof of Theorem \SUMthm
\Proof Let $f$ and $g$ be strong representors for $S$ and $T$, respectively.   Then (\VZone) is satisfied with $H = \ts\bigcupn_{\lambda > 0}\lambda\[A\(\D(S)\) - \D(T)\] = \bigcupn_{\lambda > 0}\lambda\[A\(\D(\M f)\) - \D(\M g)\]$.
Let $k$ be as in (\VZtwo).   From Theorem \VZthm(a), $\kbar$ is a strongly representative function, from which $\M \kbar$ is a strongly representable multifunction.   Since $\kbar^* = k^*$, (\VZthree) and the fact that $A^{\T\T}\wh x = \wh{Ax}$ imply that
$$\eqalign{\kbar^@(x,x^*) &= \kbar^*(x^*,\wh x) = k^*(x^*,\wh x) = \minn_{y^* \in F^*}\[f^*(x^* -  A^{\T}y^*,\wh x) + g^*(y^*,\wh{Ax})\]\cr
&= \minn_{y^* \in F^*}\[f^@(x,x^* -  A^{\T}y^*) + g^@(Ax,y^*)\].}$$
From three applications of Lemma \FFATlem(a), $\M f^@ = \M f = S$, $\M g^@ = \M g = T$ and $\M \kbar = \M \kbar^@$.   Thus
$$\eqalign{x^* &\in \(\M \kbar\)x \iff x^* \in \(\M \kbar^@\)x\cr
&\iff \exs\ y^* \in F^*\ \st\ x^* -  A^{\T}y^* \in (\M f^@)x\ \and\ y^* \in (\M g^@)Ax\cr
&\iff \exs\ y^* \in F^*\ \st\ x^* -  A^{\T}y^* \in Sx\ \and\ y^* \in TAx\cr
&\iff x^* \in \(S + A^{\T}TA\)x.}$$
This completes the proof of Theorem \SUMthm.\qed
\medbreak
We obtain the next result by taking $F = E$ and $A$ the identity map.
%:  Corollary \SUMcor
\defCorollary \SUMcor
\medbreak
\noindent
{\bf Corollary \SUMcor.}\enspace\slant Let $E$ be a Banach space, $S\colon E \toto E^*$ and $T\colon E \toto E^*$ be strongly representable and $\ts\bigcupn_{\lambda > 0}\lambda\[\D(S) - \D(T)\]$ be a closed linear subspace of E.  Then the multifunction $S + T$ is strongly representable.\endslant
\medskip
Before proceeding with the dual versions of Theorem \SUMthm\ and Corollary \SUMcor, we need to introduce some further notation.
%:  Notation \Dnot
\defNotation \Dnot
\medbreak
\noindent
{\bf Notation \Dnot.}\enspace If $S\colon\ E \toto E^*$ then we define $\Stilde\colon\ E^* \toto E\dst$ by the rule:
$$(x^*,x\dst) \in \G\(\Stilde\) \iff \infn_{(s,s^*) \in \G(S)}\bra{s^* - x^*} {\wh{s} - x\dst} \ge 0.$$
If $S$ is maximally monotone then clearly
$$(x^*,\wh x) \in \G\(\Stilde\) \iff (x^*,x) \in \G\(S^{-1}\).\meqno\NOTone$$
\($\Stilde$ is the inverse of the map $\Sbar$ defined originally by Gossez in \cite\GOSSEZ\endcite.\)
\par
We will use the following result proved in Marques Alves--Svaiter, \cite\MASUE, Theorem 1.4, pp.\ 412--413\endcite:
%:  Theorem \RESthm
\defTheorem \RESthm
\medskip
\noindent
{\bf Theorem \RESthm.}\enspace\slant Let $f \in \PCLSC(E \times E^*)$ be strongly representative.   Then $\M f^* = \wt{\M f}$.\endslant
%:  Theorem \PARthm
\defTheorem \PARthm
\medbreak
\noindent
{\bf Theorem \PARthm.}\enspace\slant Let $E$ and $F$ be Banach spaces, $S\colon E \toto E^*$ and $T\colon F \toto F^*$ be strongly representable, $B \in \L(F,E)$ and $\ts\bigcupn_{\lambda > 0}\lambda\[B^{\T}\(\R(S)\) - \R(T)\]$ be a closed linear subspace of $F^*$.  Then:
\par
\noindent
{\rm(a)}\enspace The multifunction $x \mapsto \(\Stilde + B^{\T\T}\Ttilde B^{\T}\)^{-1}\wh x$ is strongly representable. 
\par
\noindent
{\rm(b)}\enspace If, in addition, $\R(\Ttilde) \subset \wh F$ then the multifunction $\(S^{-1} + BT^{-1}B^{\T}\)^{-1}$ is strongly representable.\endslant\ (Of course, this additional condition is automatic if $F$ is reflexive.)
\Proof Let $f$ and $g$ be strong representors for $S$ and $T$, respectively.   Then (\VZfour) is satisfied with $H = \ts\bigcupn_{\lambda > 0}\lambda\[B^{\T}\(\R(S)\) - \R(T)\] = \bigcupn_{\lambda > 0}\lambda\[B^{\T}\(\R(\M f)\) - \R(\M g)\]$.
Let $k$ be as in (\VZfive).   Arguing as in Theorem \SUMthm, only using Theorem \VZthm(b) and (\VZsix) instead of Theorem \VZthm(a) and (\VZthree),
$$\kbar^@(x,x^*) = \minn_{y\dst \in F\dst}\[f^*(x^*,\wh x -  B^{\T\T}y\dst) + g^*(B^{\T}x^*,y\dst)\].$$
Consequently, from Lemma \FFATlem(a) and two applications of Theorem \RESthm,
\locno\PARone
$$\eqalignno{&x^* \in \(\M\kbar\)x \iff x^* \in \(\M\kbar^@\)x\cr
&\iff \exs\ y\dst \in F\dst\ \st\
\wh x -  B^{\T\T}y\dst \in (\M f^*)x^*\ \and\ y\dst \in (\M g^*)B^{\T}x^*\cr
&\iff \exs\ y\dst \in F\dst\ \st\
\wh x -  B^{\T\T}y\dst \in \Stilde x^*\ \and\ y\dst \in \Ttilde B^{\T}x^*&(\PARone)\cr
&\iff \wh x \in \(\Stilde + B^{\T\T}\Ttilde B^{\T}\)x^*.}$$
This completes the proof of (a).   In order to prove (b), we proceed as in (a) up to (\PARone).   The additional assumption that $\R(\Ttilde) \subset \wh F$, the fact that $B^{\T\T}\wh y = \wh{By}$, and two applications of (\NOTone) imply that
$$\eqalign{x^* \in \(\M\kbar\)x
&\iff \exs\ y \in F\ \st\ \wh x -  \wh{By} \in \Stilde x^*\ \and\ \wh y \in \Ttilde B^{\T}x^*\cr
&\iff \exs\ y \in F\ \st\ x -  By \in S^{-1}x^*\ \and\ y \in T^{-1}B^{\T}x^*\cr
&\iff x \in \(S^{-1} + BT^{-1}B^{\T}\)x^*.}$$
This completes the proof of (b).\qed
\medskip
We obtain the next result by taking $F = E$ and $B$ the identity map.
%:  Corollary \PARcor
\defCorollary \PARcor
\medbreak
\noindent
{\bf Corollary \PARcor.}\enspace\slant Let $E$ be a Banach space, $S\colon E \toto E^*$ and $T\colon E \toto E^*$ be strongly representable and $\ts\bigcupn_{\lambda > 0}\lambda\[\R(S) - \R(T)\]$ be a closed linear subspace of $E^*$.  Then:
\par
\noindent
{\rm(a)}\enspace The multifunction $x \mapsto \(\Stilde + \Ttilde\)^{-1}\wh x$ is strongly representable. 
\par
\noindent
{\rm(b)}\enspace If, in addition, $\R(\Ttilde) \subset \wh E$ then the multifunction $\(S^{-1} + T^{-1}\)^{-1}$ (known as the ``parallel sum of $S$ and $T$'') is strongly representable.\endslant\   (Of course, this additional condition is automatic if $E$ is reflexive.)
%:  Section \QUALsec More on the sandwiched closed subspace conditions
\defSection \QUALsec
\medskip
\par
\centerline{\bf\QUALsec.\quad More on the sandwiched closed subspace conditions}
\medskip
\noindent
In this section, we ``complete'' the proof of Theorem \VZthm.   We start off with a well known result on translating a strongly representative function.   See, for instance, \cite\MASBR, Proposition 3.2, p.\ 699\endcite.   \(There is a similar result in a much more general framework in \cite\HBM, Lemmas 19.13, 35.4, pp.\ 82--83, 141\endcite.\)
%:  Lemma \TRANSlem
\defLemma \TRANSlem
\medbreak
\noindent
{\bf Lemma \TRANSlem.}\enspace\slant
Let $E$ be a Banach space, $f \in \PCLSC(E \times E^*)$ be strongly representative and $(\xi,\xi^*) \in E \times E^*$.  Then we define $f_{(\xi,\xi^*)} \in \PCLSC(E \times E^*)$ by
$$f_{(\xi,\xi^*)}(x,x^*) := f(x + \xi,x^* + \xi^*) - \bra{x}{\xi^*} - \bra{\xi}{x^*} - \bra{\xi}{\xi^*}.$$
Then $f_{(\xi,\xi^*)}$ is strongly representative, $\dom\,f_{(\xi,\xi^*)} = \dom\,f - (\xi,\xi^*)$, $\G\(\M f_{(\xi,\xi^*)}\) = \G\(\M f\) -  (\xi,\xi^*)$ and $\G\(\M {f_{(\xi,\xi^*)}}^*\) = \G\(\M f^*\) -  (\xi^*,\wh{\xi})$.\endslant
\medskip
Our next result is purely algebraic in character.   In fact, Lemma \DIFFlem\ is equivalent to the known fact that if $C$ is convex then $a\in C$ and $b\in \hbox{icr}\,C \lr \,]a,b\,]\subset \hbox{icr}\,C$.   \(See Z\u{a}linescu, \cite\ZBOOK, p. 3\endcite.\)
%:  Lemma \DIFFlem
\defLemma \DIFFlem
\medbreak
\noindent
{\bf Lemma \DIFFlem.}\enspace\slant Let $C$ be a convex subset of a vector space $E$, and $H := \bigcup_{\lambda > 0}\lambda C$ be a linear subspace of $E$.   Let $x \in C$ and $\alpha \in \,]0,1[\,$.   Then
$$\ts\bigcupn_{\lambda > 0}\lambda\[C - \alpha x\] = H.\meqno\DIFFone$$\endslant
%:     Proof of Lemma \DIFFlem
\Proo $C - \alpha x \subset H - H = H$, which gives the inclusion ``$\subset$'' in (\DIFFone).   Now let $y \in H$. Then there exist $\mu > 0$ and $z \in C$ such that $y = \mu z$.   Thus
$$(1 - \alpha)z = \[(1 - \alpha)z + \alpha x\] - \alpha x \in C - \alpha x$$
and so
$$y = \mu z \in \f{\mu}{1 - \alpha}\[C - \alpha x\]
\subset \ts\bigcupn_{\lambda > 0}\lambda\[C - \alpha x\],$$
which gives the inclusion ``$\supset$'' in (\DIFFone), and thus completes the proof of Lemma \DIFFlem.\qed
\medskip
The hypotheses in Theorem \SANDthm(a) below (apart from the introduction of $\alpha$) are exactly those of Theorem \VZthm(a).   Theorem \SANDthm(b) uses the additional hypothesis that $\R\(\M g^*\) \subset \wh F$ from Theorem \PARthm(b).  
%:  Theorem \SANDthm
\defTheorem \SANDthm
\medbreak
\noindent
{\bf Theorem \SANDthm.}\enspace\slant Suppose that $E$ and $F$ are Banach spaces and that $f \in \PCLSC(E \times E^*)$ and $g \in \PCLSC(F \times F^*)$ are strongly representative.   Let $\alpha \in \,]0,1[\,$.
\par
\noindent
{\rm(a)}\enspace Suppose that $A \in \L(E,F)$ and there exists a closed linear subspace $H$ of $F$ such that
$$A\(\D(\M f)\) - \D(\M g) \subset H \subset \ts\bigcupn_{\lambda > 0}\lambda\[A(\pi_E\,\dom\,f)- \pi_F\,\dom\,g\].$$
Then
$$\alpha\[A(\pi_E\,\dom\,f)- \pi_F\,\dom\,g\] \subset A\(\D(\M f)\) - \D(\M g)\meqno\SANDone$$
and
$$\ts\bigcupn_{\lambda > 0}\lambda\[A\(\D(\M f)\) - \D(\M g)\] = \bigcupn_{\lambda > 0}\lambda\[A(\pi_E\,\dom\,f)- \pi_F\,\dom\,g\].\meqno\SANDtwo$$
\par\noindent
{\rm (b)}\enspace Suppose that $B \in \L(F,E)$, there exists a closed linear subspace $H$ of $F^*$ such that
$$B^{\T}\(\R(\M f)\) - \R(\M g) \subset H \subset \ts\bigcupn_{\lambda > 0}\lambda\[B^{\T}(\pi_{E^*}\,\dom\,f) - \pi_{F^*}\,\dom\,g\],$$
and $\R\(\M g^*\) \subset \wh F$.   Then
$$\alpha\[B^{\T}(\pi_{E^*}\,\dom\,f) - \pi_{F^*}\,\dom\,g\] \subset B^{\T}\(\R(\M f)\) - \R(\M g).\meqno\SANDthree$$
and
$$\ts\bigcupn_{\lambda > 0}\lambda\[B^{\T}(\pi_{E^*}\,\dom\,f) - \pi_{F^*}\,\dom\,g\] = \ts\bigcupn_{\lambda > 0}\lambda\[B^{\T}\(\R(\M f)\) - \R(\M g)\].\meqno\SANDfour$$
\endslant
%:     Proof of Theorem \SANDthm
\Proo(a)\enspace Let $z \in A(\pi_E\,\dom\,f)- \pi_F\,\dom\,g$.   Then there exist $(\xi,\xi^*) \in \dom\,f$ and $(\eta,\eta^*) \in \dom\,g$ such that $z = A\xi - \eta$.   Let $f_1 := f_{\alpha(\xi,\xi^*)}$ and $g_1 := g_{\alpha(\eta,\eta^*)}$.   From Lemma \TRANSlem, $f_1$ and $g_1$ are strongly representative.   Furthermore, the proof of Theorem \VZthm(a) implies that $\bigcupn_{\lambda > 0}\lambda\[A(\pi_E\,\dom\,f) - \pi_F\,\dom\,g\] = H$.   Since, from Lemma \TRANSlem\ again,
$$\eqalign{\ts\bigcupn_{\lambda > 0}\lambda\[A(\pi_E\,\dom\,f_1) - \pi_F\,\dom\,g_1\]
&= \ts\bigcupn_{\lambda > 0}\lambda\[A(\pi_E\,\dom\,f - \alpha\xi) - (\pi_F\,\dom\,g - \alpha\eta)\]\cr
&= \ts\bigcupn_{\lambda > 0}\lambda\[A(\pi_E\,\dom\,f) - \pi_F\,\dom\,g - \alpha z\],}$$
it follows from Lemma \DIFFlem\ that $\bigcupn_{\lambda > 0}\lambda\[A(\pi_E\,\dom\,f_1) - \pi_F\,\dom\,g_1\] = H$.   We now define
$$k_1(x,x^*) := \infn_{y^* \in F^*}\[f_1(x,x^* - A^{\T}y^*) + g_1(Ax,y^*)\].$$
Then Theorem \VZthm(a) shows that $\overline{k_1}$ is a strongly representative function, from which $\M \overline{k_1}$ is a strongly representable multifunction.   Thus Theorem \MASVZthm\ implies that $\D\(\M \overline{k_1}\) \ne \emptyset$, and the argument of Theorem \SUMthm\ provides $x \in \D(\M f_1) = \D(\M f) - \alpha\xi$ such that $Ax \in  \D(\M g_1) = \D(\M g) - \alpha\eta$.   But then
$$\alpha z = A(\alpha\xi) - \alpha\eta \in A\(\D(\M f) - x\) - \(\D(\M g) - Ax\) = A\(\D(\M f)\) - \D(\M g).$$
This gives (\SANDone), and (\SANDtwo) follows from  (\SANDone), (\VZone) and the observation that
$$\eqalign{\ts\bigcupn_{\lambda > 0}\lambda\[A(\pi_E\,\dom\,f) - \pi_F\,\dom\,g\]
&\subset \ts\bigcupn_{\lambda > 0}(\lambda/\alpha)\[A\(\D(\M f)\) - \D(\M g)\]\cr
&= \ts\bigcupn_{\lambda > 0}\lambda\[A\(\D(\M f)\) - \D(\M g)\].}$$
\par
(b)\enspace Let $z^* \in B^{\T}(\pi_{E^*}\,\dom\,f) - \pi_{F^*}\,\dom\,g$.   Then there exist $(\xi,\xi^*) \in \dom\,f$ and $(\eta,\eta^*) \in \dom\,g$ such that $z^* = B^{\T}\xi^* - \eta^*$.   Let $f_1$ and $g_1$ be as in (a).   Then the proof of Theorem \VZthm(b) implies that $\ts\bigcupn_{\lambda > 0}\lambda\[B^{\T}(\pi_{E^*}\,\dom\,f) - \pi_{F^*}\,\dom\,g\] = H$.   Arguing exactly as in (a), $\ts\bigcupn_{\lambda > 0}\lambda\[B^{\T}(\pi_{E^*}\,\dom\,f_1) - \pi_{F^*}\,\dom\,g_1\] = H$.   We now define
$$k_1(x,x^*) := \infn_{y \in F}\[f_1(x - By,x^*) + g_1(y,B^{\T}x^*)\].$$
Then Theorem \VZthm(b) shows that $\overline{k_1}$ is a strongly representative function, from which $\M \overline{k_1}$ is a strongly representable multifunction.   Thus Theorem \MASVZthm\ implies that $\R\(\M \overline{k_1}\) \ne \emptyset$.   From Theorem \RESthm\ and Lemma \TRANSlem, $\R\(\wt{\M g_1}\) = \R\(\M {g_1}^*\) = \R\(\M g^*\) - \alpha\wh{\eta} \subset \wh F$ and the argument of Theorem \PARthm(b) taken together with two applications of Lemma \TRANSlem\ provides us with $x^* \in \R(\M f_1) = \R(\M f) - \alpha\xi^*$ such that $B^{\T}x^* \in \R(\M g_1) = \R(\M g) - \alpha\eta^*$.   But then
$$\alpha z^* = B^{\T}(\alpha\xi^*) - \alpha\eta^* \in B^{\T}\(\R(\M f) - x^*\) - \(\R(\M g) - B^{\T}x^*\) = B^{\T}\(\R(\M f)\) - \R(\M g).$$   This gives (\SANDthree), and (\SANDfour) follows exactly as in (a).\qed 
\bigbreak
\centerline{\bf References}
\nmbr\ARENS
\item{[\ARENS]} R. Arens, \slant Operational calculus of linear relations\endslant, Pacific J. Math. {\bf 11} (1961) 9--23.
\nmbr\GOSSEZ
\item{[\GOSSEZ]}J.- P. Gossez, \slant Op\'erateurs monotones non lin\'eaires dans les espaces de Banach non r\'eflexifs\endslant, J. Math. Anal. Appl. {\bf 34} (1971), 371--395.
\nmbr\MASBR
\item{[\MASBR]}M. Marques Alves and B. F. Svaiter, \slant Br\o ndsted--Rockafellar property and maximality of monotone operators representable by convex functions in non--reflexive Banach spaces\endslant, J. of Convex Anal., {\bf 15} (2008), 693--706. 
\nmbr\MASUE
\item{[\MASUE]}-----, \slant Maximal monotone operators with
a unique extension to the bidual\endslant, J. of Convex Anal.,  {\bf 16} (2009), 409--421. 
\nmbr\MASNI
\item{[\MASNI]}-----, \slant A new old class of maximal monotone operators\endslant, J. of Convex Anal.,  {\bf 16} (2009), 881--990. 
\nmbr\MASNR
\item{[\MASNR]}-----, \slant Maximal Monotonicity, Conjugation
and the Duality Product in Non--Reflexive Banach Spaces\endslant, J. of Convex Anal.,  {\bf 17} (2010), 553--563.
\nmbr\MASD
\item{[\MASD]}-----, \slant On Gossez type (D) maximal monotone operators.\endslant, J. of Convex Anal.,  {\bf 17} (2010), 1077--1088.
\nmbr\RCA
\item{[\RCA]}R. T. Rockafellar, \slant Convex analysis\endslant, Princeton Mathematical Series, {\bf 28} (1970) Princeton University Press, Princeton, N.J.
\nmbr\RS
\item{[\RS]}B. Rodrigues and S. Simons, \slant Conjugate functions and subdifferentials in nonnormed situations for operators with complete graphs\endslant, Nonlinear Anal. {\bf 12} (1988), 1069Ð-1078.
\nmbr\PANDM
\item{[\PANDM]}S. Simons, \slant Positive sets and Monotone sets\endslant, J. of Convex Anal., {\bf 14} (2007), 297--317
\nmbr\HBL
\item{[\HBL]}-----, \slant The Hahn--Banach--Lagrange theorem\endslant, Optimization, {\bf 56} (2007), 149--169.
\nmbr\HBM
\item{[\HBM]}-----, \slant From Hahn--Banach to monotonicity\endslant, 
Lecture Notes in Mathematics, {\bf 1693},\break second edition, (2008), Springer--Verlag.
\nmbr\QUADARCHIV
\item{[\QUADARCHIV]}-----, \slant Quadrivariate versions of the AttouchÐ-Brezis theorem and strong representability\endslant, arXiv:0809.0325v1, posted September 1, 2008.  
\nmbr\SSDMON
\item{[\SSDMON]}-----, \slant Banach SSD Spaces and Classes of Monotone Sets\endslant, J. of Convex Anal., {\bf 18} (2011), 227--258.
\nmbr\SZNZ
\item{[\SZNZ]}S. Simons and C. Z\u{a}linescu, \slant Fenchel duality,
Fitzpatrick functions and maximal monotonicity\endslant, J. of
Nonlinear and  Convex Anal., {\bf 6} (2005), 1--22.
\nmbr\VZJCA
\item{[\VZJCA]}M. D. Voisei and C. Z\u{a}linescu, \slant Strongly--representable monotone operators\endslant, J. of Convex Anal.,  {\bf 16} (2009), 1011--1033. 
\nmbr\ZBOOK
\item{[\ZBOOK]} C. Z\u{a}linescu, \slant Convex analysis in
general vector spaces\endslant, (2002), World Scientific.
%Sends the final list of names to the log file
\Signoff
\bigskip
{\parindent0pt
Department of Mathematics\par
University of California\par
Santa Barbara\par
CA 93106, USA\par
E-mail: simons@math.ucsb.edu}
\bye